\documentclass[12pt]{elsarticle}

\usepackage[utf8]{inputenc}
\usepackage{times}
\usepackage{mathptmx}
\usepackage{amsfonts}

\usepackage{amsmath}
\usepackage{amsthm}

\usepackage{graphicx}

\usepackage{hyperref}
\usepackage{xspace}

\newtheorem{lemma}      {Lemma}[section]

\newtheorem{theorem}    {Theorem}[section]

\newcommand{\abs}[1]    {\left\lvert#1 \right\rvert}  

\DeclareMathOperator{\lcm}{lcm}

\newcommand{\NN}{\mathbb{Z}_+}
\newcommand{\GG}{U}
\newcommand{\VV}{V}

\newcommand{\defi}[1]{\emph{#1}}

\newcommand{\CMJ}{\textrm{FCMJ}\xspace}
\newcommand{\JRP}{\textrm{JRP}\xspace}

\begin{document}

\title { {\LARGE \bfseries On the complexity of a maintenance problem for hierarchical systems} }

\author[add2]{Andreas S. Schulz}
\ead{andreas.s.schulz@tum.de}

\author[add1]{Claudio Telha~\corref{cor1}}
\ead{catelha@miuandes.cl}

\address[add2]{Department of Mathematics and School of Management, Technische Universität München, München, Germany}
\address[add1]{Facultad de Ingeniería y Ciencias Aplicadas, Universidad de los Andes, Santiago de Chile, Chile}

\cortext[cor1]{Corresponding author}

\date  { \today }

\begin{abstract}
We prove that a maintenance problem on frequency-constrained maintenance jobs with a hierarchical structure is integer-factorization hard. This result holds even on simple systems with just two components to maintain. As a corollary, we provide a first hardness result for Levi et al.'s modular maintenance scheduling problem~\cite{Levi_Magnanti_Muckstadt_Segev_Zarybnisky_2014}.

\end{abstract}

\maketitle

\newcommand{\comp}{\mathcal{C}}






The Frequency-Constrained Maintenance Jobs problem, or \CMJ problem,  is a discrete optimization problem that arises in scheduling preventive maintenance jobs for multi-component systems featuring a hierarchical assembling structure and  multiple setup costs. 
Examples of such systems include aircraft, power generators, and large production machines~\cite{Levi_Magnanti_Muckstadt_Segev_Zarybnisky_2014}. The \CMJ was introduced in~\cite{van_Dijkhuizen_2000, vanDijkhuizenvanHarten1997}. Later, it was significantly extended to submodular setup cost structures~\cite{Levi_Magnanti_Muckstadt_Segev_Zarybnisky_2014}.

\CMJ exhibits the following trade-off: to save on setup costs related to the disassembling of the system, one can group components together to be maintained simultaneously. However, by doing so, some components will be maintained more often than they should, which can be expensive. Identifying the optimal trade-off seems challenging, and solution approaches for the \CMJ have focussed on heuristics.

The computational complexity of the \CMJ problem is not known. In this note, we prove that it is integer-factorization hard. Therefore, \CMJ is unlikely to admit a polynomial-time algorithm. The concept of integer-factorization hardness was recently used in~\cite{SchulzTelha2022} to derive a hardness result for the joint replenishment problem. Regarding the hierarchical structure, our hardness result is as good as it can be: we prove that solving \CMJ on systems with two components and one module is already integer-factorization hard.



Section~\ref{sec:model} defines the \CMJ model and the extension of this problem given in~\cite{Levi_Magnanti_Muckstadt_Segev_Zarybnisky_2014}.  In Section~\ref{sec:IF}, we recall the notion of integer-factorization hardness. Section~\ref{sec:proof} contains an overview of the hardness proof, followed by the full details.

\newcommand{\nocomp}{\mathcal{M}}
\newcommand{\rootnode}{r}

\section{The \CMJ model}\label{sec:model}


The \CMJ model can be seen as a maintenance model on a hierarchical system. The formulation presented here is a slight variation of the formulations presented in~\cite[Ch.~5.1]{Levi_Magnanti_Muckstadt_Segev_Zarybnisky_2014} and~\cite{van_Dijkhuizen_2000, vanDijkhuizenvanHarten1997}. We discuss the differences between these variations at the end of this section.

We consider a system formed by a set of components $c \in \comp$ having a hierarchical structure.  The hierarchical structure is defined by a directed tree $T=(\comp \cup \nocomp,\mathcal{E})$ with edges oriented away from the root node (an arborescence). The components in $\comp$ are the leaves of the tree. The nodes in $\nocomp$ are called \defi{modules}. One module $\rootnode \in \nocomp$ is the root of the tree and represents the entire system. An arc from $v \in \nocomp$ to $v' \in \comp \cup \nocomp$ indicates that module $v$ contains component or module $v'$. The path in the tree from the root $\rootnode$ to a component $c$ describes the order in which modules must be removed to maintain $c$.

The components are to be periodically maintained on a rolling horizon basis. At time $0$, all the components are in a ``fresh state.'' Each component $c$ will be maintained, that is, set back to the fresh state, at every integer multiple of a cycle time $q_c\in \NN$. The cycle time $q_c$ is to be determined and cannot be too large: an upper bound $f_c$ on the cycle time ensures that the maintenance is effective enough for the component. We call $f_c$ the cycle limit.

Let $C_m$ be the set of components that are children of module $m\in \nocomp$ in the tree~$T$. The maintenance of $c\in C_m$ triggers additional work on $m$, such as disassembling of the module $m$. For simplicity, we call this work a maintenance as well. It is assumed that the simultaneous maintenance of one or more components in $C_m$ will always result in a single maintenance of $m$ (for example, when disassembling module $m$ is necessary to perform maintenance on a subset of components in $C_m$).

\begin{figure}[h]
\centering
\includegraphics[scale=1.3]{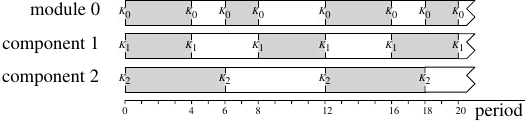}
\caption[\CMJ\ model]{An example of the \CMJ\ model. One module $0$ contains two components $1$ and $2$ with cycle limits $f_1=5$ and $f_2=6$. The diagram depicts a feasible schedule where components have cycle times $q_1=4$ and $q_2=6$.  For a given module or component  $j$, the alternating gray and white areas represent time intervals between maintenances. The setup cost $K_j$ is incurred whenever $j$ is maintained.}\label{cmj1}
\end{figure}

There are two types of costs associated with the maintenance of these components. One type is component-based.  A component $c$ incurs a setup cost $K_c$ per maintenance, thus leading to a cost rate of  $K_c/q_c$ per period. The second type is module-based. A module $m\in \nocomp$ incurs a setup cost $K_m$ each period we maintain $m$ (i.e., we maintain one or more components that are children of $m$). However, the cost rate of a module has a more complex expression because the maintenance schedule of a module does not necessarily have a constant cycle time. See the maintenance schedule of module 0 in Figure~\ref{cmj1} for an example. In general, a module $m\in \nocomp$ has a cost rate
$$K_m \sum_{\emptyset \neq C \subseteq C_m}\frac{ (-1)^{\abs{C}+1}}{\lcm\left(\left\{q_c: c\in C \right \}\right) },  $$
where $\lcm(S)$ is the least common multiple of a set of  numbers $S$. One can derive this expression for the cost rate using the inclusion-exclusion principle.

The objective of \CMJ\ is to determine the maintenance period $q_c$ of each component $c\in \comp$ to minimize the setup cost rate of all components and modules while respecting the cycle limit $f_c$. That is,

\begin{align*}
\min \ & \sum_{c \in \comp}{\frac{K_c}{q_c}} + \sum_{m \in \nocomp}K_m\sum_{C \subseteq C_m}\frac{ (-1)^{\abs{C}+1}}{\lcm(q_c: c\in C) }\\
\text{subject to}\quad  &q_c \leq f_c \qquad \text{for all } c \in \comp \tag{\CMJ},\\
  &q_c  \in\NN \qquad \text{for all } c \in \comp .
 \end{align*}

The \CMJ problem has appeared in two papers with subtle, but significant differences. In~\cite{vanDijkhuizenvanHarten1997}, the components with identical cycle times are considered part of the same cluster. A module $m$ is maintained once for each cluster that has a component in $C_m$ being maintained. This reduces the problem to find the optimal clustering. In \cite{Levi_Magnanti_Muckstadt_Segev_Zarybnisky_2014}, there are only components (no modules). At a given period, the maintenance cost of a set of components $S$ is a submodular function on $S$. The cost function of \CMJ corresponds to the special case of a submodular function defined by a dependency tree with path-additive costs~\cite[Ch~5.1]{Levi_Magnanti_Muckstadt_Segev_Zarybnisky_2014}. We mention these papers in particular because they impose cycle limits for components. The joint overhaul problem and other problems with economies of scale for grouping components are also related to the \CMJ problem. See~\cite{Dekker_Wildeman_Van_Der_Duyn_Schouten_1997} and~\cite{Nicolai_Dekker_2008} for more references.

\section{Integer-factorization hardness}\label{sec:IF}
Integer factorization is the problem of decomposing a composite number into a product of prime numbers. Finding this prime factorization can be computationally challenging, but the problem is not known to be NP-hard. The construction of several cryptographic systems relies on the difficulty of integer factorization to keep the encrypted information secure. Examples include RSA~\citep{Rivest_Shamir_Adleman_1978} and other cryptographic systems~\citep{Blum_Goldwasser_1985,Rabin_1979,SchmidtSamoa_2006}.

The notion of integer-factorization hardness for optimization problems was introduced in~\citep{SchulzTelha2022}.  An optimization problem $\Pi$ is integer-factorization hard if there is a polynomial-time algorithm that computes the factorization of a composite number $M$ using one or more oracle calls to an algorithm that solves instances of $\Pi$, counting oracle calls as $O(1)$.  If there is a polynomial-time algorithm for $\Pi$, then there is a polynomial-time algorithm for integer factorization, which is considered unlikely.

In~\citep{SchulzTelha2022}, the authors show that two variants of the joint replenishment problem are integer-factorization hard. Here we provide a similar statement, but for the \CMJ problem:

\begin{theorem}\label{thm:main}
  Suppose there is a polynomial-time exact algorithm for \CMJ.  Then, there is a polynomial-time algorithm for integer factorization.
\end{theorem}

The following section provides the proof of Theorem~\ref{thm:main}.

\section{Integer-factorization hardness of \CMJ}\label{sec:proof}


 Before proceeding, let us summarize the main ideas of the proof. We will show that  a special case of {\CMJ} is integer-factorization hard. This special case is a simple hierarchy where one module, labeled 0, contains two components labeled 1 and 2. The cost of module 0 is $K_0 \left( 1/q_1 + 1/q_2 - 1/\lcm(\{ q_1 , q_2\}) \right)$. The formulation of this special case is given in \eqref{eq:specialcase}:

\begin{align}
\min \ & \frac{K_0+K_1}{q_1} + \frac{K_0+K_2}{q_2} - \frac{K_0}{\lcm(q_1,q_2)}\notag\\
\text{subject to }  & q_1 \leq f_1, \tag{SC} \label{eq:specialcase}\\
& q_2 \leq f_2, \notag \\
& q_1, q_2 \in\NN . \notag
 \end{align}


The general idea of the reduction is the following.  Finding the integer factorization of a composite number reduces to finding one of its non-trivial divisors~\citep[Chapter~10.6, pp. 222]{Arora_Barak_2009}. Given a composite number $M$, we build an  instance~\eqref{eq:cmj3} of the maintenance problem to find one such divisor. This instance can be built in polynomial-time (polynomial in $\log(M)$) and is of the form~$\eqref{eq:specialcase}$.


Let us denote by $(q_1,q_2)$ an arbitrary optimal solution of \eqref{eq:cmj3}. The argument behind the reduction is based on a two-step process, similar to that used in~\cite{SchulzTelha2022} for the joint-replenishment problem. First, we show that $q_1= M$. Then, we show $\gcd(M,q_2)$ is a non-trivial divisor of $M$.

\subsection{Step 1}

We only need to set some of the parameters for this step. We set $K_0=1$ to normalize the coefficients, and we set $K_2=0$ as the fixed cost of this component is not needed. We also assign a high maintenance cost $K_1=M^2(M-1)-1$ to component 1, which has a cycle limit of $f_1=M$. We aim to ensure that component 1 has a cycle time equal to $M$. The model becomes:

\begin{align}
\min \ & \frac{1+K_1}{q_1} + \frac{1}{q_2} - \frac{1}{\lcm(q_1,q_2)}\notag\\
\text{subject to }  & q_1 \leq M,  \label{eq:step1}\\
& q_2 \leq f_2, \notag\\
& q_1, q_2 \in\NN . \notag
 \end{align}

Let $\GG(q_1,q_2)$ be the objective function of \eqref{eq:step1}. Note that for any $q_1\leq  M-1$ we have $$\GG(q_1,q_2)\geq (1+K_1)/q_1 \geq M^2.$$

On the other hand, for $q_1= M$ we have
$$\GG(q_1,q_2)\leq   (1+K_1)/q_1 + 1 \leq M^2 - M + 1.$$

Therefore, $\GG(M,q_2) < \GG(q_1,q_2)$ for any $q_1\leq  M-1$. It follows that $q_1$ must be equal to $M$ in any optimal solution to~\eqref{eq:step1}.

\subsection{Step 2}

As a result of Step 1, we can set $q_1=M$ and ignore the now constant terms in~\eqref{eq:step1}. We set $f_2=M-2$. We obtain the following optimization problem:
\begin{align}\label{eq:sys777}
\min\phantom{dd} & \frac{1}{q_2} - \frac{1}{\lcm(M,q_2)}\notag\\
\text{subject to \phantom{dd}} &  q_2\leq M-2, \\
\text{\phantom{s.t\phantom{dd} }} & q_2 \in\NN . \notag
 \end{align}


If $M$ is composite, we will show that any optimal solution $q_2$ of \eqref{eq:sys777} is not coprime with $M$.

\begin{lemma}\label{hardness:lemma24}
Suppose that $M$ is a composite number. If $q_2$ is coprime with $M$, then $q_2$ is not an optimal solution to (\ref{eq:sys777}).
\end{lemma}
\begin{proof}

Let $\VV(q_2)$ be the objective function of~\eqref{eq:sys777}. If $q_2$ is coprime with $M$, then $\lcm(M,q_2)=M q_2$ and, therefore,

\begin{equation*}\label{hard98}
\VV(q_2) \geq  \min_{q\leq f_2}\left\{ \frac{1}{q} - \frac{1}{Mq}\right\}= \left( \frac{M-1}{M}\right) \frac{1}{M-2} > \frac{1}{M}.
\end{equation*}

To show that any such solution $q_2$ is not optimal, we prove the existence of another feasible solution whose objective function value is smaller than $1/M$. Because $M$ is composite, we can write  $M=ab$, where $a,b\geq 2$. Let $q_2\equiv a(b-1)=M-a$. Clearly, $q_2$ is a feasible solution to (\ref{eq:sys777}), and we will prove that $\VV(q_2)<1/M$. Since $\lcm(M,q_2)\leq M(b-1)$, it follows that
\begin{equation*}
\VV(q_2)=\frac{1}{q_2} - \frac{1}{\lcm(M,q_2)} \leq \frac{1}{q_2} - \frac{1}{M(b-1)} = \frac{1}{M-a} - \frac{1}{M(b-1)} = \frac{1}{M},
\end{equation*}
as desired.

\end{proof}

Therefore, for an optimal solution $q_2$ of \eqref{eq:sys777}, $\gcd(q_2, M)$ will be greater than~1. In addition, $\gcd(q_2, M)$ will be smaller than $M$ because $q_2\leq M-2$. It follows that $\gcd(q_2, M)$ will be a non-trivial divisor of $M$.

\subsection{Finalizing the reduction}
It is well known that the greatest common divisor of two numbers can be computed in polynomial time.  If the \CMJ\ problem can be solved in polynomial time, then we can find a non-trivial divisor of a composite number $M$ in polynomial time by solving the following instance of the \CMJ problem,
\begin{align}
\min \ & \frac{M^2(M-1)}{q_1} + \frac{1}{q_2} - \frac{1}{\lcm(q_1,q_2)}\notag\\
\text{subject to }  & q_1 \leq M \tag{$\CMJ_M$}, \label{eq:cmj3}\\
& q_2 \leq M-2, \notag\\
& q_1, q_2 \in\NN, \notag
 \end{align}
and then computing $\gcd(q_2, M)$ for the optimal value of $q_2$. We have shown the following result.

\begin{theorem}
The \CMJ\ problem is integer-factorization hard, even on simple hierarchical structures with only two components in one module.
\end{theorem}

\section{Conclusions}
Note that the hierarchical structure (one module containing two components) cannot be simplified any further while keeping integer-factorization hardness. In simpler hierarchical systems, the problem is trivial.

The hardness result for the \CMJ\ problem and the hardness result for the aperiodic joint replenishment problem (\JRP) (given in~\cite[Thm. 1(b)]{SchulzTelha2022}) share a common proof structure and hold for small instances that have several commonalities. One may wonder whether there is a direct reduction from the small instances of the aperiodic \JRP that are known to be integer factorization hard (see~\cite[Lemma~1]{SchulzTelha2022} for the description of such instances) to \CMJ. However, finding a mathematical relation between optimal solutions to small \JRP instances and small instances of \CMJ is not trivial. In particular, this work suggests that the frequency constraint of Eq.~\eqref{eq:sys777} and the holding cost of~\cite[Eq.5]{SchulzTelha2022} should play a corresponding role in the mathematical relation (since in the \JRP and the \CMJ reductions they act as a soft and hard upper bound for the only variable of the model, respectively). This approach would seem to require substantially more work than the proof presented here.

\section{Acknowledgments}
 C. Telha gratefully acknowledges the support of ANID FONDECYT 11200616 and Programa Regional STIC AMSUD 22-STIC-09.

\bibliographystyle{plain}

\end{document}